\documentclass[draftcls,12pt,onecolumn]{IEEEtran}

\usepackage{times}
\usepackage{amsmath}
\usepackage{amsfonts}
\usepackage{amssymb}
\usepackage[pdftex]{graphicx}
\DeclareGraphicsExtensions{.pdf,.jpeg,.png}
\usepackage{cite}
\usepackage{url}
\usepackage{pdfpages}
\usepackage{epstopdf}
\usepackage{caption}
\usepackage{subcaption}

\title{\LARGE \bf Local Weak Observability Conditions of Sensorless AC Drives}
\author{Mohamad Koteich, \textit{Student Member, IEEE}, \\ Abdelmalek Maloum, Gilles Duc and Guillaume Sandou
	\thanks{Mohamad Koteich is with Renault S.A.S. Technocentre, 78288 Guyancourt, France, and also with L2S - CentraleSup\'{e}lec - CNRS - Paris-Sud University, 91192 Gif-sur-Yvette, France (e-mail: mohamad.koteich@renault.com).}
	\thanks{Abdelmalek Maloum is with Renault S.A.S. Technocentre, 78288 Guyancourt, France (e-mail: abdelmalek.maloum@renault.com).}
	\thanks{Gilles Duc and Guillaume Sandou are with L2S - CentraleSup\'{e}lec - CNRS - Paris-Sud University, 91192 Gif-sur-Yvette, France  (e-mail: gilles.duc@centralesupelec.fr; guillaume.sandou@centralesupelec.fr).}}

\begin{document}
\maketitle
\thispagestyle{empty}
\pagestyle{empty}
\pagenumbering{arabic}
\section*{Keywords}
\begin{center}
Sensorless control, Synchronous motor, Induction Motor, AC machine.
\end{center}

\section*{Abstract}

Alternating current (AC) electrical drive control without mechanical sensors is an active research topic. This paper studies the observability of both induction machine and synchronous machine sensorless drives. Observer-based sensorless techniques are known for their deteriorated performance in some operating conditions. An observability analysis of the machines helps understanding (and improving) the observer's behavior in the aforementioned conditions.

\section{Introduction}
Eco-friendly technologies have attracted worldwide attention due to several environmental issues. In this context, electric motors have become a serious competitor with combustion engines in many industrial applications, especially in automotive industry.

\medskip
The use of AC drives has been pioneered by the recent advances in power semiconductor switching frequencies, and power converter topologies \cite{bose}. High performance control of AC drives can be achieved using vector control \cite{leonhard}, which requires the measurement of the rotor position. For many reasons, mainly for cost reduction and reliability increase, mechanical sensorless techniques have attracted the attention of researchers as well as many large manufacturers \cite{vas}. These techniques consist of sensing the motor currents and voltages, and using them as inputs to an estimation algorithm (such as the state-observer algorithm) that estimates the rotor angular speed and/or position.

\medskip
One limitation of the use of sensorless techniques is the deteriorated performance in some operating conditions: namely the low-speed operation in the case of synchronous machines (SMs), and the low-frequency input voltage in the case of induction machines (IMs) \cite{holtz}. Usually, this problem is viewed as a stability problem, and sometimes is treated based on experimental results. However, the real problem lies in the so-called ``observability conditions'' of the system. 

\medskip
Over the past few years, a promising approach, based on the local weak observability concept \cite{hermann}, has been used in order to better understand the deteriorated performance of sensorless AC drives. Several papers have been published on this topic: authors in  \cite{canudas} \cite{ghanes} \cite{marino} \cite{vaclavek} \cite{glumineau} study the local observability of IMs in some operating conditions, namely for zero acceleration and low-frequency input voltages. Among SMs, only the permanent magnet synchronous machine (PMSM) is studied in the literature   \cite{vaclavek} \cite{glumineau} \cite{zhu} \cite{ezzat10} \cite{zaltni10}.

\medskip
A unified approach of AC drives observability analysis is proposed in this paper. Sufficient conditions for both IM and SM's observability are studied. Concerning the IM, a more general study, that covers a wider region of operating conditions, is done; the speed is not considered to be constant, which requires the estimation of the resistant torque. Concerning the SMs, the wound-rotor synchronous machine (WRSM) is studied, and considered to be the general case of SMs; its model and observability conditions are easily extended to the permanent magnet synchronous machine (PMSM) and the synchronous reluctance machine (SyRM).

\medskip
The main purpose of this paper is to contribute to a better understanding of the deteriorated performance of sensorless techniques. The paper is divided into five sections; after this introduction, the local weak observability concept is presented in section 2. Sections 3 and 4 are dedicated to the observability study of IMs and SMs respectively. Conclusions are made in section 5.

\medskip

\section{Local weak observability theory}
The \emph{local weak observability} concept \cite{hermann}, based on the rank criterion, is introduced in this section.

\subsection{Problem statement}
Systems of the following form (denoted $\Sigma$) are considered:
\begin{equation}
\Sigma :\left\{
\begin{aligned}
\dot{x} &= f\left(x(t), u(t)\right)\\
y &= h\left(x(t)\right)
\label{sigma}
\end{aligned}
\right.\end{equation}
where $x \in X \subset \mathbb{R}^n$ is the state vector, $u \in U \subset \mathbb{R}^m$ is the control vector (input), $y \in \mathbb{R}^p $ is the output vector, $f$ and $h$ are $C^\infty$ functions.

\medskip
The observation problem can be then formulated as follows \cite{besancon}: \emph{Given a system described by a representation \eqref{sigma}, find an accurate estimate $\hat{x}(t)$ for $x(t)$ from the knowledge of $u(\tau)$, $y(\tau)$ for $0 \leq \tau \leq t$}.

\medskip
The system observability is required for observer design. For nonlinear systems, \emph{global} observability is not practical, since the observer often requires the system to be \emph{instantaneously} observable in a certain \emph{neighborhood} of the state trajectories; the system should be \emph{locally weakly observable}.

\subsection{Observability rank condition}
The system $\Sigma$ is said to satisfy the observability rank condition at $x_0$, if the observability matrix, denoted by $\mathcal{O}_y(x)$, is full rank at $x_0$. $\mathcal{O}_y(x)$ is given by:
\begin{equation}
\mathcal{O}_y(x) = \frac{\partial}{\partial x}\left[ \begin{matrix}
\mathcal{L}^0_fh(x) & 
\mathcal{L}_fh(x) & 
\mathcal{L}_f^2h(x) & 
\ldots & 
\mathcal{L}_f^{n-1}h(x)
\end{matrix} \right]_{x=x_0}^T
\end{equation}
where $\mathcal{L}_f^{k}h(x)$ is the $k$th-order \emph{Lie derivative} of the function $h$ with respect to the vector field $f$. 

\subsection{Observability theorem }
A system $\Sigma$ \eqref{sigma} satisfying the observability rank condition at $x_0$ is locally weakly observable at $x_0$. More generally, a system $\Sigma$ \eqref{sigma} satisfying the observability rank condition for any $x_0$, is locally weakly observable \cite{hermann}. Rank criterion gives only a sufficient condition for local weak observability.

\section{Observability analysis of induction machines}
This section deals with the local observability conditions of induction machines. First, the machine model is presented, then its local weak observability is studied. In this study, the speed is not considered to be constant, therefore the resistant torque, which is assumed to vary slowly, should be added to the estimated state vector.

\subsection{Machine model} 
\begin{figure}[!ht]
 \centering
 \begin{subfigure}{.4\textwidth}
   \centering
   \includegraphics[width=\linewidth]{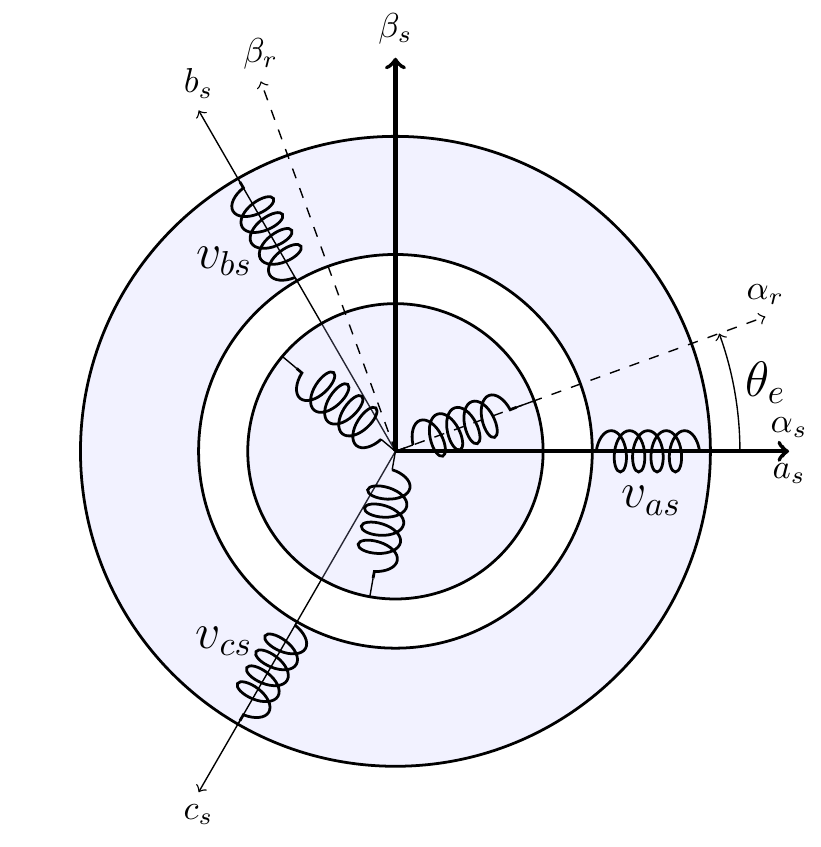}
   \caption{Schematic representation}
   \label{im_schema}
 \end{subfigure}%
 \begin{subfigure}{.6\textwidth}
   \centering
   \includegraphics[width=\linewidth]{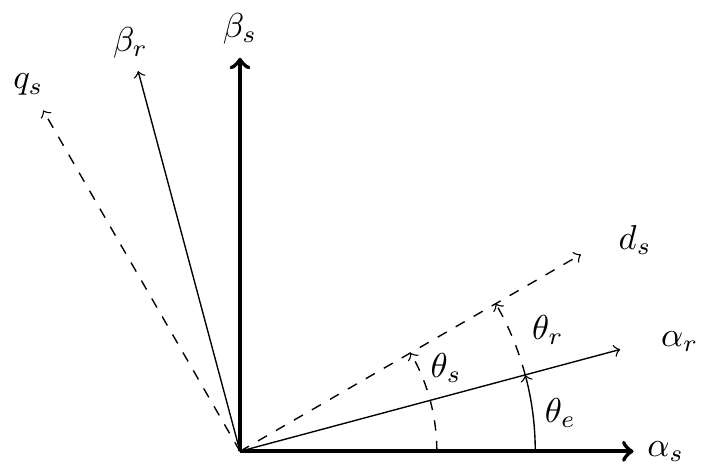}
   \caption{Vector diagram}
   \label{im_vector}
 \end{subfigure}
 \caption{Schematic representation (a) and vector diagram (b) of the induction machine. ($\alpha_s \beta_s$), ($\alpha_r \beta_r$) and ($d_s q_s$) are respectively the two-phase stator, rotor, and rotating magnetic field reference frames.}
 \label{im}
 \end{figure}

The IM  can be modeled as a three-phase stator and a three-phase rotor. The stator is supplied by a three-phase source, whereas the rotor windings are in short-circuit (see Figure \ref{im}). The IM state-space model in the two-phase stator reference frame ($\alpha_s \beta_s$) can be written as follows:
\begin{subequations}
\begin{eqnarray}
\frac{d \mathcal{I}_s}{dt} &=& \frac{1}{\sigma L_s} \left[\mathcal{V}_s - r_s \mathcal{I}_s  + \frac{M}{L_r} \left(\frac{1}{\tau_r} \mathbb{I}_2 - \omega_e \mathbb{J}_2\right) {{\Psi}_r} \right] \label{model_i}\\
\frac{d{\Psi}_r}{dt} &=& - \left(\frac{1}{\tau_r} \mathbb{I}_2 - \omega_e \mathbb{J}_2\right) {\Psi}_r + \frac{M}{\tau_r} \mathcal{I}_s \label{model_psi}\\
\frac{d \omega_e}{dt} &=& \frac{3}{2} \frac{p^2}{J} \frac{M}{L_r}  \mathcal{I}_s^T\mathbb{J}_2 \Psi_r - \frac{p}{J} T_r \label{model_w}\\
\frac{d T_r}{dt} &=& 0
\label{model_Tr}
\end{eqnarray}
\label{ss_im}
\end{subequations}
where
\begin{equation}
r_s = R_s + R_r \frac{M^2}{L_r^2} ~~~;~~~ \sigma = 1 - \frac{M^2}{L_s L_r} ~~~;~~~ \tau_r = \frac{L_r}{R_r}
\end{equation}
$R$, $L$ and $M$ stand respectively for the resistance, inductance and mutual inductance. The indices $s$ and $r$ stand for stator and rotor parameters. $p$ is the number of pole pairs, $J$ is the inertia of the rotor with the associated load, $\omega_e$ is the electrical speed of the rotor and $T_r$ is the resistant torque. ${\mathcal{I}}_s$, ${\Psi}_r$ and $\mathcal{V}_s$ stand for the stator currents, rotor fluxes and stator voltages in the stator reference frame ($\alpha_s \beta_s$) :
\begin{equation}
\mathcal{I}_s = \begin{bmatrix}
{i}_{s\alpha_s} & {i}_{s\beta_s}
\end{bmatrix}^T~~~~;~~~~
{\Psi}_r  = \begin{bmatrix}
{\psi}_{r\alpha_s} & {\psi}_{r\beta_s}
\end{bmatrix}^T~~~~;~~~~
{\mathcal{V}}_s = \begin{bmatrix}
{v}_{s\alpha_s} & {v}_{s\beta_s} 
\end{bmatrix}^T
\end{equation}
The model \eqref{ss_im} can be fitted to the structure $\Sigma$ \eqref{sigma} by taking:

\begin{eqnarray}
x = \begin{bmatrix}
{\mathcal{I}}_s^T & {\Psi}_r^T & \omega_e & T_r
\end{bmatrix}^T
~~~~;~~~~
u = \mathcal{V}_s
~~~~;~~~~
y = {\mathcal{I}}_s
\label{ss_vectors}
\end{eqnarray}

\medskip
$\mathbb{I}_n$ is the $n \times n$ identity matrix, and $\mathbb{J}_2$ is the $\pi/2$ rotation matrix :
\begin{eqnarray}
\mathbb{I}_2 = \begin{bmatrix}
1 & 0\\ 0 & 1
\end{bmatrix}~~~~~~~~~;~~~~~~~~~\mathbb{J}_2 = \begin{bmatrix}
0 & -1 \\ 1 & 0
\end{bmatrix}
\end{eqnarray}
\subsection{Observability study} 
The system \eqref{ss_im}  is a 6-${th}$ order system. The local observability study requires the evaluation of derivatives up to the 5-${th}$ order of the output. However, regarding the equations complexity, only the first and second order derivatives are evaluated for the IM in this paper.
Higher order derivatives of the output are too lengthy, and very difficult to deal with. As the observability rank criterion provides sufficient conditions, information contained in the first and second order derivatives is rich enough to study the IM observability. It should be noticed that, in practice, observer-based techniques do not usually need information on higher than second order derivatives. 

\medskip
Symbolic math software is used te evaluate complex expressions. Moreover, the following change of variables is made in order to make the study easier:
\begin{eqnarray}
\widetilde{\mathcal{I}}_s &=& \sigma L_s \mathcal{I}_s\\
\widetilde{\Psi}_r &=& \frac{M}{L_r} \Psi_r
\end{eqnarray}

The system \eqref{ss_im} becomes:
\begin{subequations}
\begin{eqnarray}
\frac{d \mathcal{\widetilde{I}}_s}{dt} &=& \mathcal{V}_s + a \mathcal{\widetilde{I}}_s  + \gamma(t) {\widetilde{\Psi}_r} \label{di_s}\\
\frac{d \widetilde{\Psi}_r}{dt} &=& - \gamma(t) {\widetilde{\Psi}_r} - \left(a-b\right) \mathcal{\widetilde{I}}_s \label{dpsi_r}\\
\frac{d \omega_e}{dt} &=& \frac{c}{J}\mathcal{\widetilde{I}}_s^T\mathbb{J}_2 \widetilde{\Psi}_r - \frac{p}{J}T_r\\
\frac{d T_r}{dt} &=& 0
\end{eqnarray}
\end{subequations}
with
\begin{eqnarray}
a = -\frac{r_s}{\sigma L_s} ~;~ b = -\frac{R_s}{\sigma L_s} ~;~ c = \frac{3 p^2}{2 \sigma L_s} ~;~ \gamma(t) = \left(\frac{1}{\tau_r} \mathbb{I}_2 - \omega_e \mathbb{J}_2\right) ~;~ \frac{d \gamma}{dt} = - \frac{d \omega_e}{dt} \mathbb{J}_2
\end{eqnarray}

\medskip
The scaled output is:
\begin{equation}
y = \mathcal{\widetilde{I}}_s = \left[
\begin{matrix}
\widetilde{i}_{s \alpha_s}\\
\widetilde{i}_{s \beta_s}
\end{matrix}\right] 
\end{equation}

its first order derivative is:
\begin{eqnarray}
\dot{y} = \frac{d\widetilde{\mathcal{I}}_s}{dt} = \mathcal{V}_s + a \mathcal{\widetilde{I}}_s  + \gamma(t) {\widetilde{\Psi}_r}
\end{eqnarray}

Adding \eqref{di_s} and \eqref{dpsi_r} gives:
\begin{eqnarray}
\frac{d\mathcal{\widetilde{I}}_s}{dt} + \frac{d\widetilde{\Psi}_r}{dt} = \mathcal{V}_s + b \mathcal{\widetilde{I}}_s
\end{eqnarray}
then:
\begin{eqnarray}
\frac{d\widetilde{\Psi}_r}{dt} = \mathcal{V}_s + b \mathcal{\widetilde{I}}_s  - \frac{d\mathcal{\widetilde{I}}_s}{dt}
\end{eqnarray}

\medskip

The second order derivative of the output can be then written as:
\begin{eqnarray}
\frac{d^2\mathcal{\widetilde{I}}_s}{dt^2} &=& \frac{d\mathcal{\widetilde{V}}_s}{dt} + a \frac{d\mathcal{\widetilde{I}}_s}{dt} + \gamma(t) \frac{d\widetilde{\Psi}_r}{dt} + \frac{d\gamma}{dt}\widetilde{\Psi}_r\\
&=& \frac{d\mathcal{V}_s}{dt} + a \frac{d\mathcal{\widetilde{I}}_s}{dt} + \gamma(t) \mathcal{V}_s + \gamma(t) b \mathcal{\widetilde{I}}_s - \gamma(t) \frac{d\mathcal{\widetilde{I}}_s}{dt} + \frac{d\gamma}{dt} \widetilde{\Psi}_r \\
&=& \frac{d\mathcal{V}_s}{dt} + \gamma(t) \mathcal{V}_s + (a \mathbb{I}_2 - \gamma(t)) \frac{d \widetilde{\mathcal{I}}_s}{dt} + \gamma(t) b \mathcal{\widetilde{I}}_s + \frac{d\gamma}{dt} \widetilde{\Psi}_r
\end{eqnarray}

\medskip
The observability study is done using the scaled output and its derivatives:
\begin{eqnarray}
y &=& \mathcal{\widetilde{I}}_s\\
\dot{y} &=& \mathcal{V}_s + a \mathcal{\widetilde{I}}_s  + \gamma(t) {\widetilde{\Psi}_r}\\
\ddot{y} &=& \frac{d\mathcal{V}_s}{dt} + a \mathcal{V}_s + \left(a^2 \mathbb{I}_2 - (a-b) \gamma(t) \right) \widetilde{\mathcal{I}}_s + \left(\frac{d\gamma}{dt} + a \gamma(t) - \gamma(t)^2\right) \widetilde{\Psi}_r 
\end{eqnarray}

\medskip
The IM observability matrix, evaluated for these derivatives, can be written as:
\begin{eqnarray}
\mathcal{O}_y^{IM} &=& \begin{bmatrix} 
1 & 0 & 0 & 0 & 0 & 0\\
0 & 1 & 0 & 0 & 0 & 0\\
a & 0 & \frac{1}{\tau_r} & \omega_e &  \widetilde{\psi}_{r \beta_s} & 0\\
0 & a & -\omega_e & \frac{1}{\tau_r} & -\widetilde{\psi}_{r \alpha_s} & 0\\  
d_{11} & d_{12} & e_{11} & e_{12} & f_{11} & f_{12}\\
d_{21} & d_{22} & e_{21} & e_{22} & f_{21} & f_{22} 
\end{bmatrix}
\label{IM_obsv_matrix}
\end{eqnarray}
with:
\begin{eqnarray}
d_{11} &=& a^2 - \frac{a-b}{\tau_r} - \frac{c}{J} \widetilde{\psi}_{r \beta_s}^2 ~~~~~~;~~~~~~ d_{12} = -(a-b) \omega_e + \frac{c}{J} \widetilde{\psi}_{r \alpha_s} \widetilde{\psi}_{r \beta_s}\\
d_{22} &=& a^2 - \frac{a-b}{\tau_r} - \frac{c}{J} \widetilde{\psi}_{r \alpha_s}^2 ~~~~~~;~~~~~~ d_{21} = (a-b) \omega_e + \frac{c}{J} \widetilde{\psi}_{r \alpha_s} \widetilde{\psi}_{r \beta_s} 
\end{eqnarray}
\begin{eqnarray}
e_{11} &=& \frac{a}{\tau_r} - \frac{1}{\tau_r^2} + \omega_e^2 + \frac{c}{J} \widetilde{i}_{s \beta_s} \widetilde{\psi}_{r \beta_s} ~~~;~~~ e_{12} = a \omega_e - 2 \frac{\omega_e}{\tau_r} + \frac{d\omega_e}{dt} - \frac{c}{J} \widetilde{i}_{s \alpha_s} \widetilde{\psi}_{r \beta_s}\\
e_{22} &=& \frac{a}{\tau_r} - \frac{1}{\tau_r^2} + \omega_e^2 + \frac{c}{J} \widetilde{i}_{s \alpha_s} \widetilde{\psi}_{r \alpha_s} ~~~;~~~ e_{21} =  - a \omega_e + 2 \frac{\omega_e}{\tau_r} - \frac{d\omega_e}{dt} - \frac{c}{J} \widetilde{i}_{s \beta_s} \widetilde{\psi}_{r \alpha_s}
\end{eqnarray}
\begin{eqnarray}
f_{11} &=& 2 \omega_e \widetilde{\psi}_{r \alpha_s} - (a-b) \widetilde{i}_{s \beta_s} + \left(a - \frac{2}{\tau_r}\right) \widetilde{\psi}_{r \beta_s}  ~~~~~~;~~~~~~
f_{12} = -\frac{p}{J} \widetilde{\psi}_{r \beta_s}\\
f_{21} &=& 2 \omega_e \widetilde{\psi}_{r \beta_s} + (a-b) \widetilde{i}_{s \alpha_s} - \left(a + \frac{2}{\tau_r}\right) \widetilde{\psi}_{r \alpha_s} ~~~~~~;~~~~~~ f_{22} = \frac{p}{J} \widetilde{\psi}_{r \alpha_s}
\end{eqnarray}

\medskip
The matrix $\mathcal{O}_y^{IM}$ \eqref{IM_obsv_matrix} is a $6\times6$ matrix, its determinant is the following:

\begin{eqnarray}
\Delta_{IM} = \frac{p}{J} \frac{M^2}{L_r^2} \left[ \frac{1}{\tau_r} \frac{d\omega_e}{dt}\left({\psi}_{r \alpha_s}^2 + {\psi}_{r \beta_s}^2\right) - \left(\omega_e^2 + \frac{1}{\tau_r^2}\right) \left( \frac{d {\psi}_{r \alpha_s}}{dt} {\psi}_{r \beta_s} - \frac{d {\psi}_{r \beta_s}}{dt} {\psi}_{r \alpha_s}\right)\right]
\label{delta_im}
\end{eqnarray}
The expression \eqref{delta_im} gives a sufficient local weak observability condition of the system \eqref{ss_im}; the value of $\Delta_{IM}$ has to be non-zero to ensure the local observability of sensorless IM drives. 

\subsubsection{Results interpretation} 
The expression \eqref{delta_im} evokes in the following two particular situations (to ensure the local weak observability of the system):
\begin{itemize}
\item $\dot{\omega}_e = 0$: For constant rotor speed, the rotor fluxes should be neither constant (or slowly varying) nor linearly dependent, i.e. $\psi_{r\alpha_s} \neq k.\psi_{r\beta_s}$, for all constant k.
\item $\dot{\psi}_{r \alpha_s} = \dot{\psi}_{r \beta_s} =0$: If the rotor fluxes are slowly varying (constant), which usually corresponds to low-frequency input voltages, the speed has to vary (rapidly).
\end{itemize}  

\subsubsection{Geometrical interpretation} The IM is locally weakly observable if the determinant $\Delta_{IM}$ \eqref{delta_im} is different from zero, this means:
\begin{eqnarray}
\frac{\tau_r \dot{\omega}_e}{1 + \tau_r^2 \omega_e^2} \neq \frac{\frac{d {\psi}_{r \alpha_s}}{dt} {\psi}_{r \beta_s} - \frac{d {\psi}_{r \beta_s}}{dt} {\psi}_{r \alpha_s}}{{\psi}_{r \alpha_s}^2 + {\psi}_{r \beta_s}^2} \Longleftrightarrow 
\frac{d}{dt} \arctan(\tau_r \omega_e) \neq \frac{d}{dt}\arctan\left(\frac{\psi_{r \beta_s}}{\psi_{r \alpha_s}}\right) \label{cond_im_0}
\end{eqnarray}
Knowing that (Figure \ref{im_vector}):
\begin{eqnarray}
\frac{d}{dt}\arctan\left(\frac{\psi_{r \beta_s}}{\psi_{r \alpha_s}}\right) = \frac{d}{dt}\theta_s
\end{eqnarray}
the condition \eqref{cond_im_0} becomes:
\begin{eqnarray}
\omega_s \neq \frac{d}{dt} \arctan(\tau_r \omega_e) \label{im_geo}
\end{eqnarray}
where $\omega_s$ is the angular frequency of the rotating magnetic field. The relationship \eqref{im_geo} tells that if the rotor speed $\omega_e$ is slowly varying, the input voltage frequency should not be very low. Nevertheless, in a wide range of applications, the rotor speed is (almost) constant during the most of machine operation time, the main conclusion to ensure the IM local weak observability is that the input voltage frequency should not be very low.

\section{Observability analysis of synchronous machines}
\begin{figure}[!ht]
\centering
\begin{subfigure}{.4\textwidth}
  \centering
  \includegraphics[width=\linewidth]{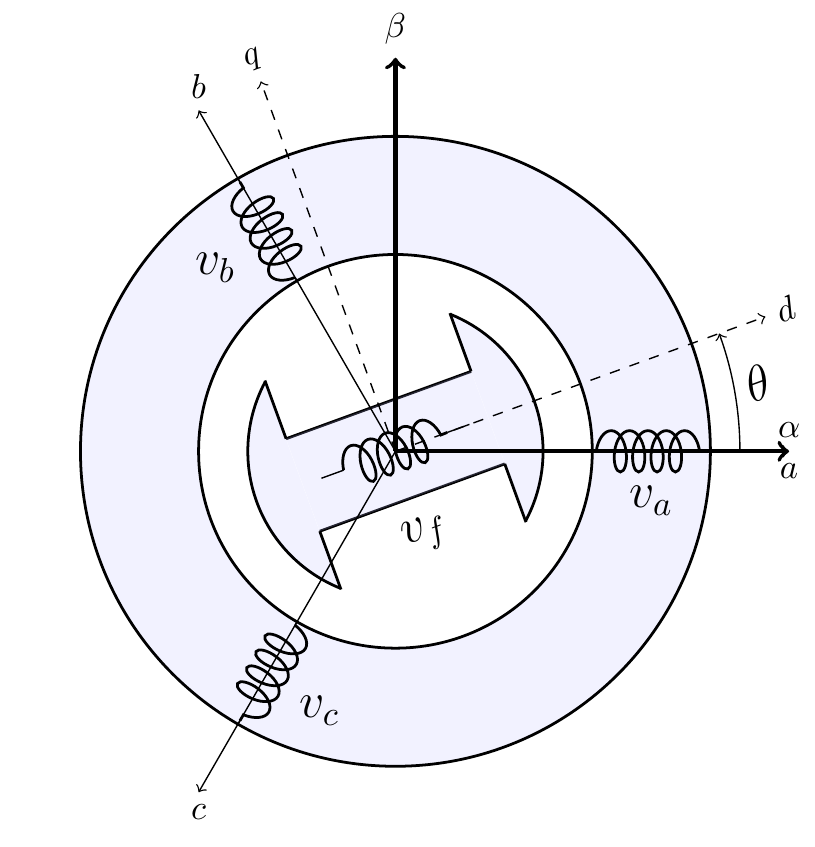}
  \caption{Schematic representation}
  \label{sm_schema}
\end{subfigure}%
\begin{subfigure}{.6\textwidth}
  \centering
  \includegraphics[width=\linewidth]{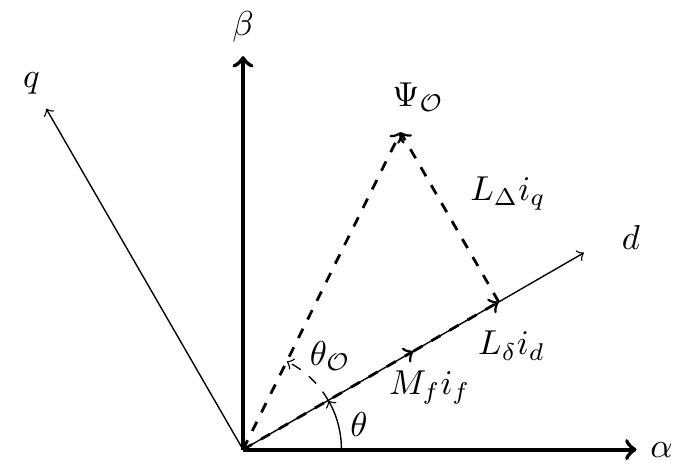}
  \caption{Vector diagram showing the (dashed) observability vector}
  \label{sm_vector}
\end{subfigure}
\caption{Schematic representation (a) and vector diagram (b) of the wound-rotor synchronous machine. ($\alpha \beta$) and ($dq$) are respectively the stator and rotor reference frames.}
\label{sm}
\end{figure}

\begin{figure}[!ht]
\centering
\begin{subfigure}{.4\textwidth}
  \centering
  \includegraphics[width=\linewidth]{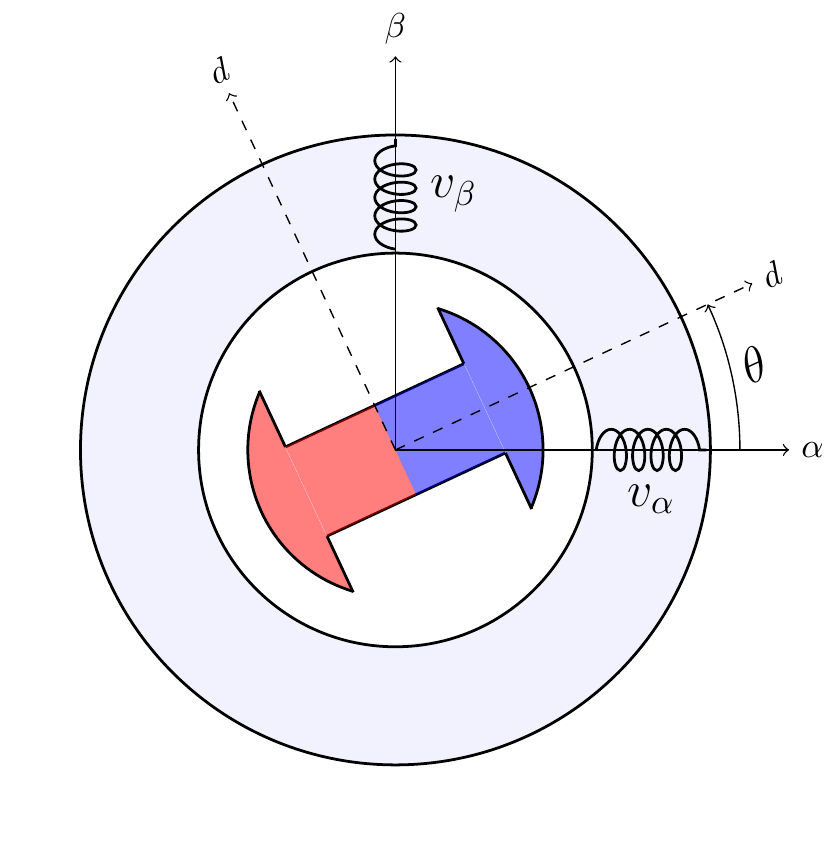}
  \caption{IPMSM}
  \label{ipmsm_schema}
\end{subfigure}%
\begin{subfigure}{.4\textwidth}
  \centering
  \includegraphics[width=\linewidth]{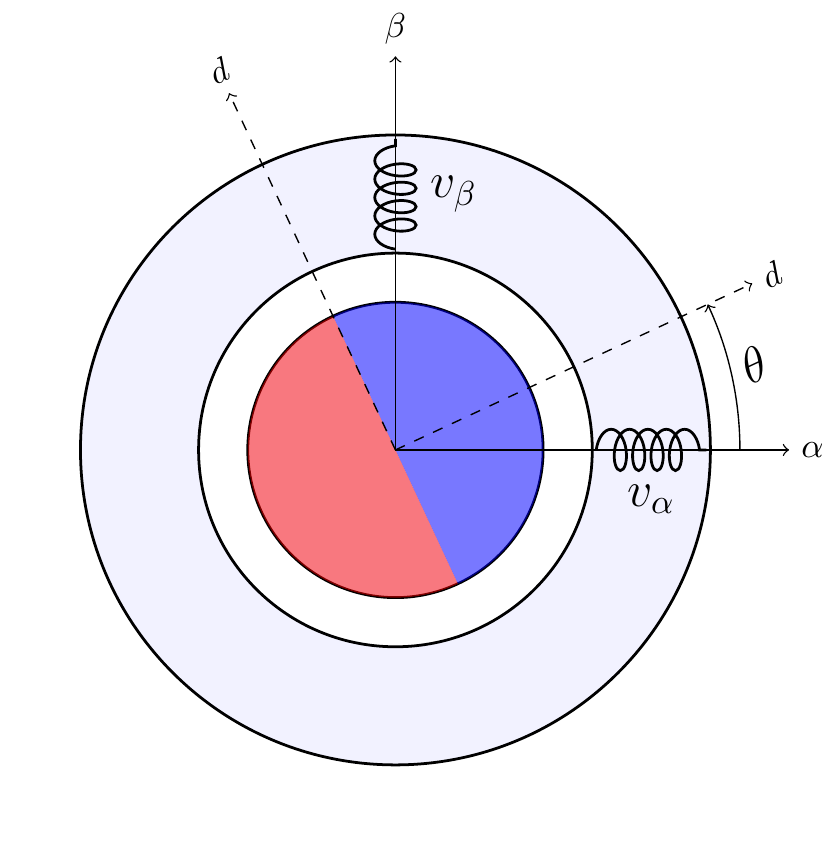}
  \caption{SPMSM}
  \label{spmsm_schema}
\end{subfigure}
\caption{Schematic representation of the IPMSM (a) and the SPMSM (b)}
\label{pmsm}
\end{figure}

\begin{figure}[!ht]
\centering
\begin{subfigure}{.4\textwidth}
  \centering
  \includegraphics[width=\linewidth]{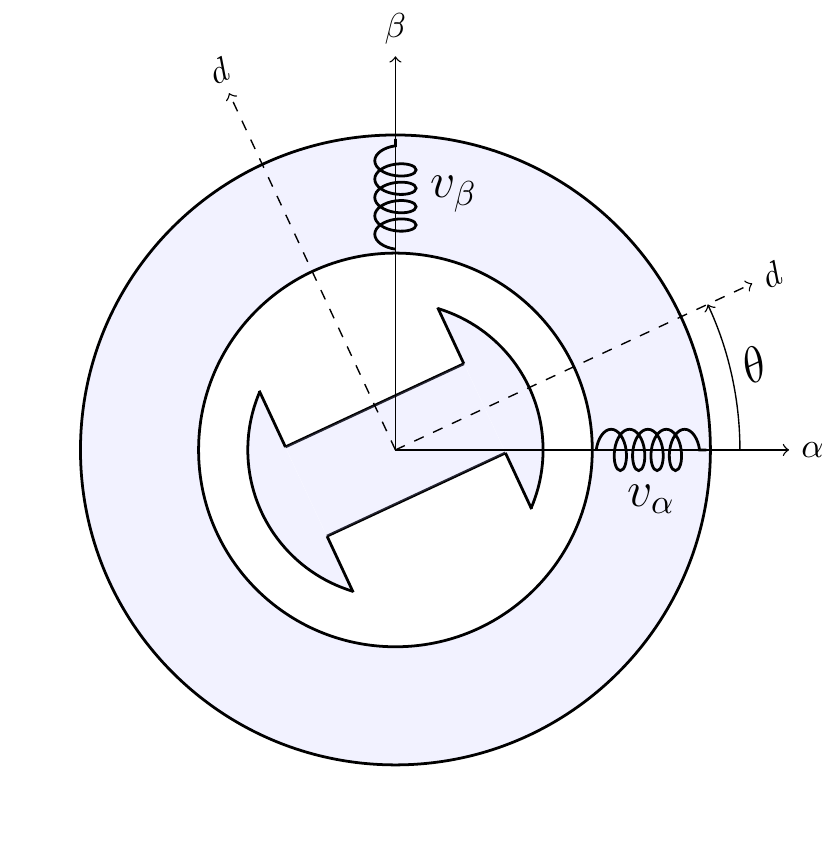}
  \caption{SyRM}
  \label{syrm_schema}
\end{subfigure}%
\begin{subfigure}{.4\textwidth}
  \centering
  \includegraphics[width=\linewidth]{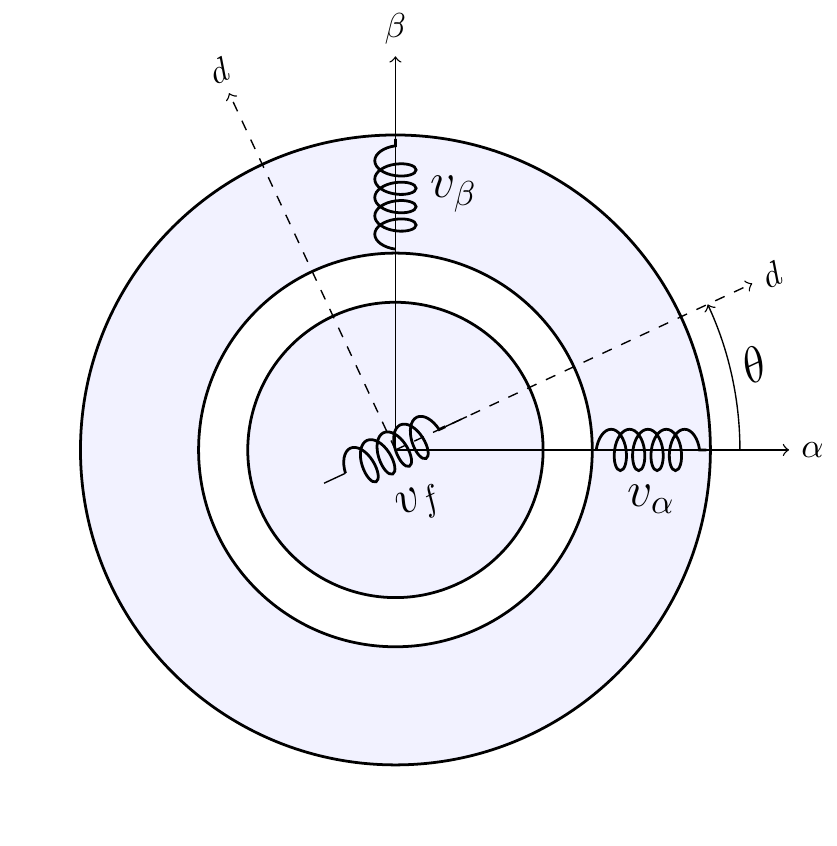}
  \caption{Non-salient WRSM}
  \label{nswrsm_schema}
\end{subfigure}
\caption{Schematic representation of the SyRM (a) and the non-salient WRSM (b)}
\label{syrm}
\end{figure}

An SM is a three-phase machine, having the same stator structure as the IM. Depending on the rotor structure, there exist three main SM sub-families (Figures \ref{sm}, \ref{pmsm}, and \ref{syrm}): wound-rotor (WRSM), permanent-magnet (PMSM) and reluctance type SyRM synchronous machines. Both WRSM and PMSM can be either salient type (non cylindrical) rotor or non-salient type (cylindrical) rotor. Interior PMSM (IPMSM) is a salient machine, surface-mounted PMSM (SPMSM) is a non-salient machine. In this section, WRSM is studied, and considered to be the general case of SMs.
\subsection{Wound-rotor synchronous machine model} The state-space model of the WRSM is written in the two-phase ($\alpha \beta$) stationary reference frame  (see Figure \ref{sm_schema}), in a way to be fitted to the structure \eqref{sigma}:
\begin{subequations}
\begin{eqnarray}
\frac{d{\mathcal{I}}}{dt}&=&-{\mathfrak{L}^{-1}}{\mathfrak{R}_{eq}}\mathcal{I}
+ {\mathfrak{L}^{-1}}\mathcal{V}\\
\frac{d\omega}{dt} &=& \frac{p}{J} T_m - \frac{p}{J} T_r\\
\frac{d\theta}{dt} &=& \omega
\end{eqnarray}
\label{ss_sm}
\end{subequations}
where the state, input and output vectors are respectively:
\begin{eqnarray}
x = \left[
\begin{matrix}
\mathcal{I}^T & \omega & \theta
\end{matrix}
\right]^T~~~;~~~
u = 
\mathcal{V}= \left[\begin{matrix}
v_\alpha & v_\beta & v_f
\end{matrix}\right]^T~~~;~~~
y = 
\mathcal{I} = \left[\begin{matrix}
i_\alpha & 
i_\beta & 
i_f
\end{matrix}\right]^T
\end{eqnarray}
$\mathcal{I}$ and $\mathcal{V}$ are the current and voltage vectors.
Indices $\alpha$ and $\beta$ stand for stator signals, index $f$ stands for rotor (\textit{field}) ones. $\omega$ stands for the electrical rotor angular speed and $\theta$ for the electrical angular position of the rotor. $\mathfrak{L}$ is the (position-dependent) matrix of inductances, and $\mathfrak{R}$ is the matrix of resistances:
\begin{equation}
\mathfrak{L}=\left[ \begin{matrix}
   {{L}_{0}}+{{L}_{2}}\cos 2\theta  & {{L}_{2}}\sin 2\theta  & {{M}_{f}}\cos \theta   \\
   {{L}_{2}}\sin 2\theta  & {{L}_{0}}-{{L}_{2}}\cos 2\theta  & {{M}_{f}}\sin \theta   \\
   {{M}_{f}}\cos \theta  & {{M}_{f}}\sin \theta  & {{L}_{f}}  \\
\end{matrix} \right] ~;~
\mathfrak{R}_{eq} = \mathfrak{R} + \frac{\partial\mathfrak{L}}{\partial\theta}\omega ~;~
\mathfrak{R} = \left[\begin{matrix}
R_s & 0   &  0 \\
  0 & R_s &  0 \\
  0 & 0   & R_f
\end{matrix}\right] \nonumber
\end{equation}
$J$ is the moment of inertia of the rotor with its associated load, $p$ is the number of pole pairs, $T_r$ is the resistant torque and $T_m$ is the motor torque:
\begin{eqnarray}
T_m = \frac{3}{2}~p M_f i_f (i_{\beta} \cos \theta - i_{\alpha} \sin \theta)  -\frac{3}{2}~pL_2 \left[(i_{\alpha}^2 - i_{\beta}^2 ) \sin 2\theta - 2 i_{\alpha} i_{\beta} \cos 2\theta \right]
\end{eqnarray}

\subsection{Other SMs models}
The other SMs can be seen as special cases of the salient-type WRSM; the IPMSM (Figure \ref{ipmsm_schema}) model can be derived by considering the rotor magnetic flux to be constant:
\begin{eqnarray}
\frac{di_f}{dt} = 0
\label{cond_pmsm_1}
\end{eqnarray}
and by substituting $M_f i_f$ by the permanent magnet flux $\psi_r$:
\begin{eqnarray}
i_f = \frac{\psi_r}{M_f}
\label{cond_pmsm_2}
\end{eqnarray}

The SyRM (Figure \ref{syrm_schema}) model can be derived from the IPMSM model by considering the rotor magnetic flux $\psi_r$ to be zero: 
\begin{eqnarray}
\psi_r \equiv 0
\label{cond_syrm}
\end{eqnarray}

The equations of the non-salient WRSM (Figure \ref{nswrsm_schema}) and SPMSM (Figure \ref{spmsm_schema}) are the same as the salient WRSM and IPMSM respectively, except that the stator self-inductances are constant and independent of the rotor position, that is:
\begin{eqnarray}
L_2 = 0 ~~~\implies~~~ L_d = L_q = L_0
\label{non_salient}
\end{eqnarray}

\subsection{Observability study}  The system \eqref{ss_sm} is a 5-${th}$ order system. Its observability study requires the evaluation of the output derivatives up to the 4-${th}$ order. In this study, only the first order derivatives are calculated, higher order ones are very difficult to calculate and to deal with. This gives the following ``partial'' observability matrix:
\begin{equation}
\mathcal{O}_{y}^{SM} = \left[ \begin{matrix}
\mathbb{I}_{3} &\mathbb{O}_{3 \times 1} &\mathbb{O}_{3 \times 1}  \\
-\mathfrak{L}^{-1}\mathfrak{R}_{eq} &-\mathfrak{L}^{-1}\mathfrak{L}' \mathcal{I} &\mathfrak{L^{-1}}' \frac{d \mathcal{I}}{dt} - \mathfrak{L^{-1}}\mathfrak{L}''\omega \mathcal{I} 
\end{matrix} \right] 
\label{obsv_matrix}
\end{equation}
where $\mathbb{O}_{n \times m}$ is an $n \times m$ zero matrix.
$\mathfrak{L}'$ and $\mathfrak{L}''$ denote, respectively, the first and second partial derivatives of $\mathfrak{L}$ with respect to $\theta$.

\medskip

It is sufficient to have five linearly independent lines of matrix \eqref{obsv_matrix} to insure the local observability of the system. The first five lines, which come from the first derivatives of $i_{\alpha}$ and $i_{\beta}$, are studied. This choice is motivated by the fact that these currents are available for measurement in all synchronous machines, the rotor current (that gives the sixth line of matrix \eqref{obsv_matrix}) does not exist for PMSM and SyRM. Another reason comes from the physics of the machine: $i_f$ is a DC signal, whereas both $i_{\alpha}$ and $i_{\beta}$ are AC signals, so it is more convenient for physical interpretation to take them together.

\medskip
Symbolic math software is used to calculate the determinant. Currents and fluxes are expressed in the rotor ($dq$) reference frame (using the Park transformation) in order to make the interpretation easier. The determinant is written, in terms of stator currents and fluxes, under the following general form:
\begin{eqnarray}
\Delta_{SM} =  \frac{L_\Delta}{L_D L_q} \left[\left(
\frac{1}{L_\Delta}\left(\psi_d - L_q i_d\right)^2 + L_\delta  i_q^2
\right)\omega + 
\frac{d\psi_d}{dt}i_q + \frac{d\psi_q}{dt}i_d - \left(\frac{di_d}{dt} \psi_q + \frac{di_q}{dt} \psi_d  \right)\right] 
\label{delta_sm}
\end{eqnarray}
with
\begin{eqnarray}
{L}_D = L_d - \frac{M_f^2 }{L_f}~;~ {L}_\Delta = L_\delta - \frac{M_f^2 }{L_f}~;~ L_d = L_0 + L_2;~;~ L_q = L_0 - L_2;~;~ L_\delta = L_d - L_q
\end{eqnarray}

The observability conditions for the different SMs can be deduced from \eqref{delta_sm}, by replacing the fluxes by their expressions in the rotor ($dq$) reference frame as the following:
\paragraph{WRSM} $\psi_d = L_d i_d + M_f i_f~~;~~\psi_q = L_q i_q$
\begin{eqnarray}
\Delta_{WRSM} &=&  \frac{1}{{L}_D L_q} \left[
\left(L_\delta i_d + M_f i_f \right)^2 + L_\Delta {L}_\delta i_q^2
\right]\omega \nonumber \\
&& + \frac{{L}_\Delta}{{L}_D L_q} \left[
\left(L_\delta \frac{di_d}{dt} + M_f \frac{di_f}{dt} \right) i_q  -  \left(L_\delta i_d + M_f i_f \right) \frac{di_q}{dt}
\right] \label{delta_wrsm} 
\end{eqnarray}
The determinant of the non-salient rotor WRSM ($L_\delta = 0$) can be then deduced:
\begin{equation}
\Delta_{N-WRSM} = \frac{M_f^2}{{L}_D L_q} \left[
i_f^2 \omega - \frac{M_f}{L_f} \left(i_q \frac{di_f}{dt} - i_f \frac{di_q}{dt} \right)
 \right]
 \label{delta_wrsm_nsp}
\end{equation}
\paragraph{PMSM} $\psi_d = L_d i_d + \psi_r ~~;~~ \psi_q = L_q i_q~~;~~{L}_D = L_d ~~;~~ {L}_\Delta = L_\delta$
\begin{eqnarray}
\Delta_{IPMSM} &=& \frac{1}{L_d L_q} \left[\left(L_\delta i_d + \psi_r \right)^2 + L_\delta^2 i_q^2 \right]\omega + \frac{L_\delta}{L_d L_q} \left[L_\delta \frac{di_d}{dt} i_q - \left(L_\delta i_d + \psi_r \right) \frac{di_q}{dt}
\right]
\label{delta_ipmsm}\\
\Delta_{SPMSM} &=& \frac{\psi_r^2}{L_0^2}\omega
\label{delta_spmsm}
\end{eqnarray}
\paragraph{SyRM} $\psi_d = L_d i_d ~~;~~ \psi_q = L_q i_q ~~;~~{L}_D = L_d ~~;~~ {L}_\Delta = L_\delta$
\begin{eqnarray}
\Delta_{SyRM} = \frac{L_\delta^2}{L_d L_q} \left[\left( i_d^2 + i_q^2 \right)\omega + \frac{di_d}{dt} i_q - i_d \frac{di_q}{dt}
\right]
\label{delta_syrm}
\end{eqnarray}

\subsubsection{Results interpretation} 
The expression of $\Delta_{WRSM}$ \eqref{delta_wrsm} shows that WRSM can be observable at standstill if the currents $i_d$, $i_q$ and $i_f$ are not constant at the same time. Non-salient WRSM \eqref{delta_wrsm_nsp}, IPMSM \eqref{delta_ipmsm} and SyRM \eqref{delta_syrm} can be also observable at standstill. However, the observability of SPMSM at standstill cannot be guaranteed.

\subsubsection{Geometrical interpretation} The observability condition $\Delta_{WRSM} \neq 0$ implies:
\begin{eqnarray}
\omega \neq \frac{\left(L_\delta i_d + M_f i_f \right) L_\Delta \frac{di_q}{dt} - 
\left(L_\delta \frac{di_d}{dt} + M_f \frac{di_f}{dt} \right) L_\Delta i_q}{\left(L_\delta i_d + M_f i_f \right)^2 + L_\Delta L_\delta i_q^2}
\end{eqnarray}
The above equation can be written as:
\begin{eqnarray}
\omega \neq \frac{\left(L_\delta i_d + M_f i_f \right)^2 + L_\Delta^2 i_q^2}{\left(L_\delta i_d + M_f i_f \right)^2 + L_\Delta L_\delta i_q^2} \frac{\left(L_\delta i_d + M_f i_f \right) L_\Delta \frac{di_q}{dt} - 
\left(L_\delta \frac{di_d}{dt} + M_f \frac{di_f}{dt} \right) L_\Delta i_q}{\left(L_\delta i_d + M_f i_f \right)^2 + L_\Delta^2 i_q^2} \nonumber
\end{eqnarray}
then
\begin{eqnarray}
\omega \neq \frac{\left(L_\delta i_d + M_f i_f \right)^2 + L_\Delta^2 i_q^2}{\left(L_\delta i_d + M_f i_f \right)^2 + L_\Delta L_\delta i_q^2} \frac{d}{dt} \arctan \left( \frac{L_\Delta i_q}{L_\delta i_d + M_f i_f} \right) 
\end{eqnarray}
The following approximation can be adopted\footnote{This approximation does not affect the observability conditions at standstill where $\omega = 0$ and currents are non-zero. In addition, this assumption is an equality for the PMSM ansd the SyRM.}:
\begin{eqnarray}
\frac{\left(L_\delta i_d + M_f i_f \right)^2 + L_\Delta^2 i_q^2}{\left(L_\delta i_d + M_f i_f \right)^2 + L_\Delta L_\delta i_q^2} \approx 1
\label{assum}
\end{eqnarray}

Thus, the WRSM observability condition can be formulated as:
\begin{eqnarray}
\omega \neq \frac{d}{dt} \arctan \left( \frac{L_\Delta i_q}{L_\delta i_d + M_f i_f} \right) 
\label{obsv_vect_0}
\end{eqnarray}
It can be seen that the above equation describes a vector, which will be called the \textit{observability vector} and denoted $\Psi_\mathcal{O}$ (Figure \ref{sm_vector}), that has the following components in the $dq$ reference frame:
\begin{eqnarray}
\Psi_{\mathcal{O}d} &=& L_\delta i_d + M_f i_f\\
\Psi_{\mathcal{O}q} &=& L_\Delta i_q
\end{eqnarray}
The condition \eqref{obsv_vect_0} becomes:
\begin{eqnarray}
\omega \neq \frac{d}{dt} \theta_{\mathcal{O}}
\label{obsv_vect_1}
\end{eqnarray}
where $\theta_{\mathcal{O}}$ is the phase of the vector $\Psi_\mathcal{O}$ in the $dq$ reference frame.

\medskip
Finally, the following SMs local weak observability condition can be stated: the observability of an SM is guaranteed if the rotational velocity of the observability vector with respect to the rotor is different from the electrical velocity of the rotor with respect to the stator. Therefore, at standstill, the observability vector should rotate and not be fixed. 

\section{Conclusions}
The local weak observability of motion-sensorless three-phase electrical drives has been studied through this paper. The study shows that the observability might not be guaranteed in some operating conditions: such as low-frequency rotor fluxes for the induction machines, and standstill operation for the synchronous machines.

\medskip
The main conclusion, concerning the local weak observability of induction machines, is that input voltage frequency should not be low to ensure the machine observability.

\medskip
The concept of \emph{observability vector} generalizes the observability conditions for synchronous machines: this vector should rotate whith respect to the rotor with a speed different from the electrical speed of the rotor with respect to the stator, to ensure the observability.

\medskip
Having the expressions \eqref{delta_im}, \eqref{delta_wrsm}, \eqref{delta_wrsm_nsp}, \eqref{delta_ipmsm} and \eqref{delta_syrm}, it seems to be interesting to combine high-frequency signal injection techniques with observer-based sensorless techniques. In other words, high-frequency voltages can be injected in the stator (or the rotor in the case of WRSM) windings, in order to guarantee the machine observability in the critical operating conditions, which improves the observer's behaviour. Concerning the SPMSM, the only way to insure its observability is to move the rotor (equation \eqref{delta_spmsm}), to make the speed different from zero.

\end{document}